\documentclass{owrart}

\usepackage{amssymb}

\begin{document}


\begin{talk}[Andr\'e Henriques, Carlo A. Rossi, Chenchang Zhu]{Giovanni Felder}%
{Elliptic gamma functions, triptic curves and ${{\mathit{SL}}}_3(\mathbb{Z})$}{Felder, Giovanni}

\newtheorem{thm}{Theorem}

\noindent The elliptic gamma function
\cite{GF_R}
is a function of three complex variables obeying
\[
\Gamma(z+\sigma,\tau,\sigma)=
\theta_0(z,\tau)\Gamma(z,\tau,\sigma),\;
\theta_0(z,\tau)=\prod_{j=0}^\infty (1-e^{2\pi i((j+1)\tau-z)})
(1-e^{2\pi i(j\tau+z)}).
\]
In \cite{GF_FV1} three-term
relations for $\Gamma$ involving $\mathit{ISL}_3(\mathbb{Z})\!=\!{{\mathit{SL}}}_3(\mathbb{Z})\!\ltimes\! \mathbb{Z}^3$ were
discovered, generalizing the
modular properties  of theta functions under
$\mathit{ISL}_2(\mathbb{Z})\!=\!{{\mathit{SL}}}_2(\mathbb{Z})\!\ltimes\! \mathbb{Z}^2$.

Here we summarize the results of \cite{GF_fhrz} where we show
that these identities are special cases
 of a set of three-term
relations for a family of gamma
functions $\Gamma_{a,b}$, which are interpreted geometrically as
giving a meromorphic section of a hermitian
gerbe on the {\em universal triptic curve}.
This result generalizes the fact that the theta
function $\theta_0$ is a section of a hermitian line bundle on the universal
elliptic curve.

First, we describe the gerbe by the enlarged gamma function family.
For $a, b$ linearly independent in  the set $\Lambda_{prim}$ of
primitive
vectors (namely not multiples of other vectors) in the lattice
$\Lambda=\mathbb{Z}^3$, there is a unique primitive $\gamma\in \Lambda^\vee_{prim}$
in the dual lattice such that
$\det (a, b, \cdot) = s \gamma$ for some $s>0$. For $w\in\mathbb{C}, x\in\Lambda\otimes\mathbb{C}=
\mathbb{C}^3$ for which the products converge we define
\[ \Gamma_{a,b}(w,x):= \prod_{\delta\in C_{+-}/\mathbb{Z}\gamma}
(1-e^{-2\pi i (\delta(x)-w)/\gamma(x)}) \prod_{\delta\in
C_{-+}/\mathbb{Z}\gamma} (1 - e^{2\pi i (\delta(x)-w)/\gamma(x)})^{-1} ,\]
where $C_{+-}=C_{+-}(a,b)=\{\delta\in\Lambda^\vee|\delta(a)>0, \delta(b)\leq 0\}$ and $\mathbb{Z}
\gamma$ acts on it by translation. We set similarly $C_{-+}(a,b)=C_{+-}(b,a)$. We define $\Gamma_{a,\pm a}=1$.
The function $\Gamma_{a,b}$ is meromorphic on $\mathbb{C}\times(U^+_a\cap U^+_b)$,
where $U^+_a=\{x\in \mathbb{C}^3| \mathrm{Im}(\alpha(x)\overline{\beta(x)})>0\}$
for any oriented basis $\alpha$, $\beta$ of the plane $\delta(a)=0$. For
linearly independent $a,b\in\Lambda_{prim}$,
$\Gamma_{a,b}$ is a 
finite product of ordinary elliptic gamma
functions: \begin{equation}\label{eq:gamma-ab-ordinary}
\Gamma_{a,b}(w,x)=\prod_{\delta\in F/\mathbb Z\gamma}
\Gamma\left(\frac{w+\delta(x)}{\gamma(x)},\frac{\alpha(x)}{\gamma(x)},
\frac{\beta(x)}{\gamma(x)}\right),
\end{equation}
for any $\alpha,\beta \in \Lambda^\vee$  satisfying $\alpha(b)=\beta(a)=0$
and $\alpha(a)>0$ $\beta(b)>0$, $F= \{\delta\in\Lambda^\vee| 0\leq\delta(a)<\alpha(a), 0\leq\delta(b)<\beta(b)\}$.
In particular, we recover $\Gamma(z,\tau,\sigma)=\Gamma_{a,b}(z,(\tau,\sigma,1))$ for 
the standard basis vectors $a=e_1$, $b=e_2$.

The functions $\Gamma_{a,b}$ satisfy cocycle conditions generalizing the three-term
relations of \cite{GF_FV1}:
\begin{gather}\label{eq:gamma}
\Gamma_{a,b}(w,x)\Gamma_{b,a}(w,x)=1, \quad  x\in U_a^+\cap
U_b^+, \\
\notag
\Gamma_{a,b}(w,x)\Gamma_{b,c}(w,x)\Gamma_{c,a}(w,x)=\exp\left(-\frac{\pi i}{3}
P_{a,b,c}(w,x)\right), \quad x\in U_a^+\cap U_b^+\cap
U_c^+,
\end{gather}
where  $P_{a,b,c}(w,x)\in \mathbb Q(x)[w]$ can be explicitly described in terms
of the Bernoulli polynomial $B_{3,3}$, see \cite{GF_fhrz}.
Moreover the gamma functions obey cocycle identities related to the
action of the group $\mathit{ISL}_3(\mathbb{Z})={{\mathit{SL}}}_3(\mathbb{Z})\ltimes\mathbb{Z}^3$. Fix a
{\em framing} of $\Lambda_{prim}$, namely for each $a\in\Lambda_{prim}$
a choice of oriented
basis $(\alpha_1, \alpha_2, \alpha_3)$ of $\Lambda^\vee\otimes\mathbb{R}$
such that $\alpha_1(a)=1$, $\alpha_2(a)=\alpha_3(a)=0$.
Let \[\Delta_a ((g, \mu); w, x)=\prod_{j=0}^{\mu( g^{-1}a)-1}
\theta_0
\left(
 \frac{w+j\alpha_1(x)}{\alpha_3(x)},\frac{\alpha_2(x)}{\alpha_3(x)}
\right), \]where $(g, \mu) \in \mathit{ISL}_3(\mathbb{Z})={{\mathit{SL}}}_3(\mathbb{Z})\ltimes\mathbb{Z}^3$.
 Then we have
\begin{equation}\label{eq:gamma-delta}
\frac{\Gamma_{g^{-1}a,g^{-1}b}(w+\mu(g^{-1}x),g^{-1}x)}{\Gamma_{a,b}(w,x)}=
e^{\pi i P_{a,b}((g,\mu);w,x)}\frac{\Delta_a((g,\mu);w,x)}{\Delta_b((g,\mu);w,x)},
\end{equation}
\begin{equation}\label{eq:delta}
\Delta_a(\hat g\hat h;w,x)=e^{2\pi i P_a(\hat g,\hat h;w,x)}
\Delta_a(\hat g;w,x)\Delta_{g^{-1}a}(\hat h;w+\mu(g^{-1}x),g^{-1}x),
\end{equation}
where $\hat g=(g,\mu), \hat h=(h,\nu)$ and
$P_{a,b}(\hat g;\cdot)$, $P_a(\hat g,\hat h;\cdot)$ are again in $\mathbb Q(x)[w]$.

Let $X{}$ be the dual of the tautological line bundle of $\mathbb{C}P^2$ restricted to $\mathbb{C}P^2-\mathbb{R}P^2$. Then $\mathit{ISL}_3(\mathbb{Z})$ acts on $X{}=(\mathbb{C}\times (\mathbb{C}^3-\mathbb{R}^3))/\mathbb{C}^\times$ by  $(g, \mu)\cdot [(w, x)]=[(w-\mu(x), g\cdot x)]$. Both $\Gamma_{a,b}$ and $\Delta_a$ descend to the quotient $X{}$. 
By equations \eqref{eq:gamma} \eqref{eq:gamma-delta} \eqref{eq:delta} we have:
\begin{thm}
There is an $\mathit{ISL}_3(\mathbb{Z})$-equivariant \v Cech 2-cocycle
\[(\phi_{a,b,c}, \phi_{a,b}, \phi_{a})=(e^{-\frac{2\pi i}{3!}P_{a,b,c}(w,x)},e^{-\frac{2\pi i}{2!}P_{a,b}((g,\mu);w,x)},e^{-2\pi i P_a((g,\mu),(h,\nu);w,x)}), \]  in $C^2_{\mathit{ISL}_3(\mathbb{Z})}(\mathcal{V}, \mathcal{O}^\times)$ where $\mathcal{V}$ is the equivariant covering of $X{}$ made up by $V_a=\{(w,x)|x\in U_a^+\}/\mathbb{C}^\times$, $a\in\Lambda_{prim}$. The image of $\phi$
in the equivariant \v Cech complex with
values in the sheaf $\mathcal{M}^\times$ of
invertible meromorphic
sections is the coboundary of the equivariant cochain $(\Gamma_{a,b},\Delta_a)
\in C^1_{\mathit{ISL}_3(\mathbb{Z})}(\mathcal{V},\mathcal{M}^\times)$.

\end{thm}
The {\bf gamma gerbe} $\mathcal{G}$ is the holomorphic equivariant
gerbe on $X$ corresponding to $\phi$. Equivalently, it is
a holomorphic gerbe on the stack
$\mathcal{X}=[X{}/\mathit{ISL}_3(\mathbb{Z})]$.

More geometrically, if we view gerbes over stacks as central extensions of group\-oids, then $\mathcal{G}$ is presented by a group\-oid $R\rightrightarrows U_0$ fitting in the central extension of group\-oids over $U_0$:
\[ 1\to \mathbb{C}^\times \times U_0 \to R \to U_1 \to 1,\]
where $U_0=\sqcup V_a$, $U_1=U_0\times_{X}(\mathit{ISL}_3(\mathbb{Z}) \times X) \times_{X} U_0$,
$R=\sqcup L_{a,b}\otimes L_b(g)^{-1}$ with $L_{a,b}$, $L_b(g)$, $(a,b\in\Lambda_{prim},
g\in\mathit{ISL}_3(\mathbb Z)$ the holomorphic $\mathbb{C}^\times$-bundles
with transition functions
$\phi_{a,b,c}\phi_{a,b,d}^{-1}$ (on $(V_a\cap V_b) \cap V_c \cap V_d$) and $\phi_{b,b'}(g,\cdot) \phi_{b,b''}^{-1}(g,\cdot)$ (on $V_b\cap V_{b'}\cap V_{b''}$)
respectively. Notice that $U_1=\cup W_{ g, a, g^{-1}b} $ where $W_{g, a,
g^{-1}b}=\{ (g, y)|y\in V_a, g^{-1} y \in V_{g^{-1}b}\}$. Then $\Gamma_{a,b}\Delta_b^{-1}$ provides a
meromorphic group\-oid homomorphism $U_1\to R$ hence $\Gamma$'s and $\Delta$'s can be viewed as a meromorphic section of $\mathcal{G}$.
A {\em hermitian structure} on a gerbe in this language is
simply a hermitian structure on the complex line bundle associated to the central
extension. A holomorphic gerbe with hermitian structure has a canonical
connective structure whose curvature represents its Dixmier--Douady class.

\begin{thm}
Using the notation in \eqref{eq:gamma-ab-ordinary}, there is a hermitian structure $h_{a,b} h_b^{-1}$
on
$\mathcal{G}$ with $
h_{a,b}(w,x)=\prod_{\delta\in F/\mathbb Z\gamma}h_3\left(\frac{w+\delta(x)}{\gamma(x)},\frac{\alpha(x)}{\gamma(x)},\frac{\beta(x)}{\gamma(x)}\right)
$ and
$h_a((g,\mu);w,x)=\prod_{j=0}^{\mu(g^{-1}a)-1}h_2\left(\frac{w+j\alpha_1(x)}{\alpha_3(x)},\frac{\alpha_2(x)}{\alpha_3(x)}\right)$, where
$h_n$ are defined
by Bernoulli polynomials: $h_n(z,\tau_1,\dots,\tau_{n-1})=
\exp\left(-(4\pi/n!)B_{n-1,n}(\zeta,t_1,\dots,t_{n-1})\right)$,
$\zeta=\mathrm{Im}\,z$, $t_j=\mathrm{Im}\,\tau_j$.
\end{thm}

Moreover, as with line bundles, we can construct the gamma gerbe $\mathcal{G}$ via (pseudo)-divisors.
A {\em triptic curve} $\mathcal{E}$ is a holomorphic
stack of the form $[\mathbb{C}/\iota(\mathbb{Z}^3)]$ with $\iota: \mathbb{Z}^3 \to \mathbb{C}$ a map of rank 2 over $\mathbb{R}$.
An {\em orientation} of a triptic curve $\mathcal{E}$ is given by a choice of a generator
of $H^3(\mathcal{E}, \mathbb{Z})\cong \mathbb{Z}$. Then the stack ${\mathcal T\! r}:=[(\mathbb{C}P^2-\mathbb{R}P^2)/{{\mathit{SL}}}_3(\mathbb{Z})]$ is
the moduli space of oriented triptic curves. The stack $\mathcal{X}=[X{}/\mathit{ISL}_3(\mathbb{Z})]$ is the
total space of the universal family of triptic curves over ${\mathcal T\! r}$. Given an \'etale
map $U\to \mathcal{E}$, let $Z_U=0\times_{\mathcal{E}} U$.
$Z_U$ is naturally a discrete subset of a principal oriented $\mathbb{R}$-bundle
on $U$.
A {\em pseudodivisor} on $U$ is a function $D: Z_U \to \mathbb{Z}$ such that if $\lim y_n =+\infty$
(resp. $-\infty$)
for a sequence $y_n$ in $Z_U$
with relatively compact image in $U$ then $\lim D(y_n)=1$ (resp. $0$).
The notion of positive/negative infinity is derived from the orientation
of the fibres of the $\mathbb{R}$-bundle.
We can globalize this to $\mathcal{X}$, namely for an \'etale map $U\to \mathcal{X}$, a {\em pseudodivisor} on $U$
is a function $D: {\mathcal T\! r}\times_{\mathcal{X}}U\to \mathbb{Z}$ such that for every point $q\to {\mathcal T\! r}$ with corresponding fibre
$\mathcal{E}=q\times_{{\mathcal T\! r}} \mathcal{X}$, the restriction 
to $q\times_{\mathcal T\! r} Z_U$ is a pseudodivisor on $U\times_{\mathcal{X}}\mathcal{E}$.
Then for two such $D_i$'s, the pushforward $p_*(D_1-D_2)$ is a divisor on $U$, hence can be used to twist
a line bundle $L$ to $L(p_*(D_1-D_2))$, where
$p: {\mathcal T\! r}\times_{\mathcal{X}} U \to U$. Using the categorical description of gerbes in \cite{GF_Bry93},
we then have

\begin{thm} The gamma gerbe $\mathcal{G}$ is a gerbe over $\mathcal{X}$ made up by the following
data: for $U$ with an \'etale open map $U\to \mathcal{X}$,
\begin{gather*}
Obj(\mathcal{G}_U)=\{(L, D)|\;\text{$L$ is a  line bundle on $U$ and $D$ a pseudodivisor on $U$}\}, \\
Mor(\mathcal{G}_U)\big((L_1,D_1)\to(L_2,D_2)\big)=
\Gamma^\times\left(U,\;\big(L_1^*\otimes L_2\big)\big(p_*(D_2-D_1)\big)\right),
\end{gather*}
the invertible holomorphic sections.

\end{thm}

We also have the following theorems calculating various cohomology groups and Dixmier--Douady classes of the gamma gerbe and of its restriction to a fibre.
\begin{thm} Let $\mathcal{E} = \mathbb{C}/\iota(\mathbb{Z}^r)$,
where $x_j=\iota(e_j)$, the images of the standard basis vectors,
 are assumed to be $\mathbb{Q}$-linearly independent and
to span $\mathbb{C}$ over $\mathbb{R}$. Then
\[ H^{i\leq r-2}(\mathcal{E},\mathcal{O}^\times)=\wedge^i(\mathbb{C}^r/(x_1,...,x_r)\mathbb{C})/\wedge^i(\mathbb{Z}^r), \; H^{r-1}(\mathcal{E},\mathcal{O}^\times)=\mathcal{E}\times \mathbb{Z},
\]
and $ H^{\geq r}(\mathcal{E},\mathcal{O}^\times)=0.$ In particular, for $r=3$, 
the groups classifying holomorphic and
topological gerbes on $\mathcal{E}$ are
$H^2(\mathcal{E},\mathcal{O}^\times)=\mathcal{E}\times \mathbb{Z}$ and $H^3(\mathcal{E},\mathbb{Z})=\mathbb{Z}$, respectively.
\end{thm}

\begin{thm}The Dixmier--Douady class  $c(\mathcal{G}|_\mathcal{E})$ of the
restriction of the gamma gerbe to $\mathcal{E}$ is a generator of
$H^3(\mathcal{E},\mathbb{Z})=\mathbb{Z}$.
\end{thm}

\begin{thm}$H^3(\mathcal{X}, \mathbb{Z})$ fits into the short exact sequence
\[ 0\to \mathbb{Z} \to H^3(\mathcal{X}, \mathbb{Z})/
\text{torsion} \to H^3(\mathbb{Z}^3, \mathbb{Z}) \cong \mathbb{Z} \to 0. \]
The image of the Dixmier--Douady class $c(\mathcal{G})\in H^3(\mathcal{X},\mathbb{Z}) $ of
the gamma gerbe is a generator of $H^3(\mathbb{Z}^3,\mathbb{Z})$.
\end{thm}
There should exist non-abelian versions of this story in
the context of q-deformed conformal field theory \cite{GF_FVKZB}.

\end{talk}

\end{document}